\documentclass[leqno,12pt]{article}

\usepackage{latexsym}
 \usepackage{graphicx}
\usepackage{amsmath}
\usepackage{amssymb}
\usepackage{mathtools}
\usepackage{cite}
\usepackage{color}
\usepackage{subfig}

\usepackage[margin=1.3in, top=1.3in, bottom=1.3in]{geometry}

\newcommand{\nc}{\newcommand}
\nc{\nt}{\newtheorem}
\nt{thm}{Theorem}[section]
\nt{cor}[thm]{Corollary}
\nt{prop}[thm]{Proposition}
\nt{lem}[thm]{Lemma}
\nt{defn}[thm]{Definition}
\nt{rem}[thm]{Remark}
\nt{exa}[thm]{Example}
\nt{ass}[thm]{Assumption}
\nt{alg}[thm]{Algorithm}
\nt{conj}[thm]{Conjecture}
\nc{\ip}[2]{\mbox{$\langle #1,#2 \rangle$}}
\nc{\pf}{\noindent{\bf Proof\ \ }}
\nc{\finpf}{\hfill{$\Box$}\linespace}
\nc{\linespace}{\vspace{\baselineskip} \noindent}
\nc{\R}{{\bf R}}
\nc{\E}{{\bf E}}
\nc{\W}{{\mathcal W}}
\nc{\M}{{\mathcal M}}
\nc{\Rn}{{\bf R}^n}
\nc{\bx}{\bar{x}}
\nc{\by}{\bar{y}}
\nc{\inT}{\mbox{\rm int}\,}
\nc{\cl}{\mbox{\rm cl}\,}
\nc{\gph}{\mbox{\rm gph}\,}
\nc{\argmin}{\mbox{\rm argmin}\,}

\def\tto{\;{\lower 1pt \hbox{$\rightarrow$}}\kern -12pt
           \hbox{\raise 2.8pt \hbox{$\rightarrow$}}\;}
\newenvironment{myequation}{\setcounter{equation}{\value{thm}}
   \begin{equation}}{\addtocounter{thm}{1}\end{equation}}

\nc{\bmye}{\begin{myequation}}
\nc{\emye}{\end{myequation}}

\begin{document}
\title{
The structure of conservative gradient fields
}
\author{
\and
A.S. Lewis
\thanks{ORIE, Cornell University, Ithaca, NY.
\texttt{people.orie.cornell.edu/aslewis} 
\hspace{2cm} \mbox{~}
Research supported in part by National Science Foundation Grant DMS-2006990.}
\and
Tonghua Tian
\thanks{ORIE, Cornell University, Ithaca, NY.  
\texttt{tt543@cornell.edu}}
}
\date{\today}
\maketitle

\begin{abstract}
The classical Clarke subdifferential alone is inadequate for understanding automatic differentiation in nonsmooth contexts.  Instead, we can sometimes rely on enlarged generalized gradients called ``conservative fields'', defined through the natural path-wise chain rule:  one application is the convergence analysis of gradient-based deep learning algorithms.  In the semi-algebraic case, we show that all conservative fields are in fact just Clarke subdifferentials plus normals of manifolds in underlying Whitney stratifications.
\end{abstract}
\medskip

\medskip


\section{Introduction}
Popular deep learning solvers like PyTorch \cite{paszke2017automatic} and  TensorFlow \cite{tensorflow2015-whitepaper} increasingly rely on automatic differentiation for gradient-based optimization algorithms.  Given an input point $x \in \Rn$, the solver returns a gradient-like vector $g \in \Rn$ that depends not just on the objective function $f$ itself but rather on its algorithmic representation.  At least when $f$ is smooth, we might hope that $g$ is the gradient 
$\nabla f(x)$, but even then, nonsmooth algorithmic ingredients may produce surprises.  For example, the formula
\[
f(s) ~=~ \big((-s)^+ + s \big)  - s^+ \qquad (s \in \R)
\]
(where $s^+ = \max\{0,s\}$) always outputs the value zero, and yet one implementation \cite[Appendix A.2]{bolte-pauwels-ad} of automatic differentiation in TensorFlow  outputs the derivative 
\bmye \label{bad}
g(s) = \left\{
\begin{array}{ll}
0 & (s \ne 0) \\
1 & (s=0).
\end{array}
\right.
\emye
For another recent survey of the same issue, see \cite{lee-yu-rival-yang}. 

Despite this disconcerting behavior, practitioners widely apply automatic differentiation to nonsmooth objective functions $f \colon \Rn \to \R$, as discussed in \cite{MR3511388}.  For the particular case of  stochastic subgradient descent algorithms, see \cite{MR4056927}.  Fortunately, as demonstrated by \cite{bolte-pauwels-conservative}, automatic differentiation at points $x \in \Rn$ typically does produce outputs $g(x) \in \Rn$ with gradient-like properties: time-dependent trajectories $x(\cdot)$ satisfy the {\em chain rule}
\bmye \label{chain}
\tfrac{d}{dt} f(x) ~=~ \ip{\tfrac{d}{dt}x}{g(x)} \quad \mbox{almost always},
\emye
thereby justifying the convergence of the stochastic subgradient descent method \cite{MR4056927}.

Thus motivated, Bolte and Pauwels \cite{bolte-pauwels-conservative} develop a novel and elegant notion of generalized derivative for locally Lipschitz objectives $f$ precisely around the chain rule.  They consider {\em conservative fields\/}:  closed set-valued mappings $G \colon \Rn \tto \Rn$ with the property that the chain rule (\ref{chain}) holds providing that we always select $g(x) \in G(x)$.  Being locally Lipschitz, $f$ is differentiable on a full-measure set 
\mbox{$\Omega \subset \Rn$}, and \cite{bolte-pauwels-conservative} shows that the value $G(x)$ for any conservative field must contain the set
\[
\{ \lim_r \nabla f(x_r) : x_r \to x,~ x_r \in \Omega \},
\]
and hence, if convex, also its convex hull $\partial f(x)$, the {\em Clarke subdifferential} \cite{clarke}.

An objective function with a conservative field is called {\em path differentiable}.  The theoretical existence question has a long history, surveyed in \cite{bolte-pauwels-conservative} and dating back to \cite{valadier}, but in practice, objectives are always path differentiable:  examples include smooth and convex functions, and their sums and differences, as well as semi-algebraic (or, more generally, tame  \cite{vandendries-miller}) functions.  The Clarke subdifferential is the minimal convex-valued conservative field for any path differentiable objective.  However, as the example (\ref{bad}) makes clear for automatic differentiation, we are forced to consider conservative fields larger than the Clarke subdifferential.  What do they look like in general?

In this work we focus on the most concrete case, where objective functions and conservative fields are semi-algebraic.  (The tame generalization is immediate, but we do not pursue it here.)  We prove an intuitive structural result, characterizing the conservative fields of an objective function as modest modifications of its Clarke subdifferential, arising simply by including normals to manifolds comprising Whitney stratifications of $\Rn$.  For example the conservative field (\ref{bad}) can arise from the stratification 
$\R = (-\infty,0) \cup \{0\} \cup (0,+\infty)$.  Thus, while the important idea of a conservative field is arrived at very differently from the notion of the Clarke subdifferential, in practice the two ideas are very close.

\section{Characterizing conservative fields}
Turning to the formal development, we consider set-valued operators on $\Rn$, by which we mean set-valued mappings $G \colon \Rn \tto \Rn$.  The {\em sum} of two operators $G$ and $H$ maps points $x \in \Rn$ to the sum $G(x)+H(x)$.  We call $G$ {\em closed} if its graph 
$\{(x,y) \in \Rn \times \Rn : y \in G(x)\}$ is closed, and {\em locally bounded} if every point in $\Rn$ has a neighborhood $\Omega \subset \Rn$ whose image $G(\Omega)$ is bounded.  A {\em selection} of $G$ is an operator whose graph is contained in the graph of $G$.  The following definition is from \cite[Lemma 2]{bolte-pauwels-conservative}.

\begin{defn} \label{definition}
{\rm
A {\it conservative field} for a locally Lipschitz function $f \colon \Rn \to \R$  is a closed, locally bounded, nonempty-set-valued operator $G$ on $\Rn$  such that all absolutely continuous curves
$x \colon [0,1] \to \Rn$ satisfy the following chain rule: for almost all $t \in [0,1]$,
\[ 
\tfrac{d}{dt} f\big(x(t)\big) ~=~ \ip{\tfrac{d}{dt}x(t)}{g} \quad \mbox{for all}~ g \in G\big(x(t)\big).
\]
}
\end{defn}

We next consider smooth stratifications of sets in a Euclidean space $\E$;  in all of our discussions of functions and manifolds, ``smooth'' simply means continuously differentiable.  For any smooth manifold 
$\M \subset \E$, we denote the tangent and normal spaces to $\M$ at any point $x \in \M$ by $T_\M(x)$ and $N_\M(x)$.  A finite collection $\W$ of disjoint smooth manifolds in $\E$ comprise a {\em Whitney stratification} (of their union) if the following condition holds for all manifolds $\M$ and $\M'$ in $\W$:
\[
\left.
\begin{array}{ll}
x_r \in \M, 		& y_r \in N_\M(x_r) \\
x_r \to x \in \M',	& y_r \to y
\end{array}
\right\}
\quad
\Rightarrow
\quad
\left\{
\begin{array}{l}
\M' \subset \cl\M \\
y \in N_{\M'}(x).
\end{array}
\right.
\]
In particular, we can associate any Whitney stratification $\W$ of $\Rn$ with a closed {\em normal operator} $\Phi_{\W}$ on $\Rn$ by setting $\Phi_\W = N_\M$ on each manifold $\M$ in $\W$.  We call $\W$ {\em semi-algebraic} if each $\M$ in $\W$ is semi-algebraic.

We can now state our result.

\begin{thm}[Semi-algebraic conservative fields]
Given a semi-algebraic locally Lipschitz function $f \colon \Rn \to \R$, a semi-algebraic set-valued operator is a conservative field for $f$ if and only if it is a closed, locally bounded, nonempty-valued selection of the sum of the Clarke subdifferential $\partial f$ and the normal operator for a semi-algebraic Whitney stratification of $\Rn$.
\end{thm}

\pf
Consider any semi-algebraic conservative field $G$ for $f$.  By \cite[Theorem 4]{bolte-pauwels-conservative}, there exists a Whitney stratification $\W$  of  $\Rn$ such that $f$ is smooth on each manifold $\M$ in $\W$, and at each point $x \in \M$, the Riemannian gradient $\nabla_\M f(x) \in T_\M(x)$ satisfies
\bmye \label{1}
G(x) ~\subset~ \nabla_\M f(x) + N_\M(x).
\emye
The proof in \cite{bolte-pauwels-conservative}, using stratification techniques from \cite{vandendries-miller}, makes clear that we can assume $\W$ to be semi-algebraic.  By \cite[Lemma 8]{Lewis-Clarke}, there exists a semi-algebraic Whitney stratification of the graph of $f$ that maps via the canonical projection 
$\Rn \times \R \to \Rn$ onto a semi-algebraic Whitney stratification  $\W'$  of $\Rn$  that is ``compatible'' with  $\W$:  in other words, each manifold in $\W$ is a union of manifolds in $\W'$.  

Now consider any point $x \in \Rn$.  There exist unique manifolds $\M$ in $\W$ and $\M'$ in $\W'$ containing $x$, and $\M'$ must be a submanifold of  $\M$, so 
\bmye \label{2}
N_\M(x) ~\subset~ N_{\M'}(x).   
\emye
By the definition of the Riemannian gradient, we have
\bmye \label{3}
\nabla_\M f(x) ~\in~ \nabla_{\M'} f(x) + N_{\M'}(x).
\emye
On the other hand, by \cite[Proposition 4]{Lewis-Clarke}, we have 
\[
\partial f(x) ~\subset~ \nabla_{\M'} f(x) + N_{\M'}(x).
\]
Since $f$ is locally Lipschitz, there exists a vector $g \in \partial f(x)$.  We deduce 
\bmye \label{4}
\nabla_{\M'} f(x) + N_{\M'}(x) ~=~ g + N_{\M'}(x) ~\subset~ \partial f(x) + N_{\M'}(x).
\emye
Combining the inclusions (\ref{1}), (\ref{2}), (\ref{3}), and (\ref{4}), we deduce
\[
G(x) ~\subset~ \partial f(x) + N_{\M'}(x),
\]
so $G$ is a selection of the sum $\partial f + \Phi_{\W'}$, as required.

Conversely, consider a semi-algebraic Whitney stratification $\W$ of $\Rn$.  By 
\cite[Corollary 2 and Proposition 2]{bolte-pauwels-conservative}, the subdifferential $\partial f$ is a conservative field for $f$.  On the other hand, for any radius  $r>0$,  the truncated normal operator 
defined by 
\[
\Phi^r_\W(x) ~=~ \Phi_\W(x) \cap rB \qquad (x \in \Rn),
\]
where $B \subset \Rn$ is the closed unit ball, is a conservative field for the identically zero function, by \cite[Theorem 3]{bolte-pauwels-conservative}.   Hence  $\partial f + \Phi^r_\W$ is a conservative field for  $f$, by \cite[Corollary 4]{bolte-pauwels-conservative}.  

Now consider any closed, locally bounded, nonempty-valued selection $G$ of the operator 
$\partial f + \Phi_\W$.  If $G$ is not a conservative field, then Definition~\ref{definition} fails for some absolutely continuous curve $x \colon [0,1] \to \Rn$.  The image $C = x([0,1])$ is compact, so since $G$ is locally bounded, the image $G(C)$ is bounded.  Since $f$ is locally Lipschitz, the image 
$\partial f(C)$ is also bounded.  We deduce $G(C) - \partial f(C) \subset rB$ for some radius $r>0$.  All points $x \in C$ therefore satisfy $G(x) \subset \partial f(x) + \Phi_\W^r (x)$, since, by assumption, for any vector $g \in G(x)$ there exists a subgradient $y \in \partial f(x)$ such that 
$g \in y + \Phi_\W (x)$, so in fact $g \in y + \Phi_\W^r (x)$.  But $\partial f + \Phi^r_\W$ is conservative, contradicting the failure of  Definition~\ref{definition}.
\finpf 

\bibliographystyle{plain}
\small
\parsep 0pt

\def\cprime{$'$} \def\cprime{$'$}

\bibliographystyle{plain}
\small
\parsep 0pt
\bibliography{adrian}

\end{document}